\documentclass[12pt]{amsart}
 \usepackage{amsmath,amsthm} 
 \usepackage{amssymb,amsfonts,amsxtra}
 \usepackage{psfrag,graphicx}
 \theoremstyle{plain}
\newtheorem{theorem}[subsection]{Theorem}
\newtheorem{prop}[subsection]{Proposition}

\theoremstyle{definition}
\newtheorem{ex}[subsection]{Example}
\newtheorem{rk}[subsection]{Remark}

\newcommand{\co}{\colon\thinspace}
\newcommand{\MCG}{\mathcal{MCG}}
\newcommand{\Z}{\mathbb{Z}}
\newcommand{\R}{\mathbb{R}}
\newcommand{\free}{\mathbb{F}_{k}}

\DeclareMathOperator{\Aut}{Aut}
\DeclareMathOperator{\Out}{Out}
\DeclareMathOperator{\GL}{GL}

\begin{document}
\title{Morse theory and conjugacy classes of finite subgroups}

\author[N.~Brady]{Noel Brady$^1$}
\address{Dept.\ of Mathematics\\
        University of Oklahoma\\
	Norman, OK 73019}
\email{nbrady@math.ou.edu}

\author[M.~Clay]{Matt Clay}
\address{Dept.\ of Mathematics\\
        University of Oklahoma\\
	Norman, OK 73019}
\email{mclay@math.ou.edu}

\author[P.~Dani]{Pallavi Dani}
\address{Dept.\ of Mathematics\\
        University of Oklahoma\\
	Norman, OK 73019}
\email{pdani@math.ou.edu}

\date{\today}

\begin{abstract}
	We construct a CAT($0$) group containing a finitely presented 
	subgroup with infinitely many conjugacy classes of finite-order elements. 
	Unlike previous examples (which were based on right-angled Artin groups) 
	our ambient CAT($0$) group does not contain any rank $3$ free abelian subgroups.  
	
	We also construct examples of groups of type $\mathrm F_n$ inside mapping class groups, Aut($\free$), and Out($\free$) which have infinitely many conjugacy classes of  finite-order elements. 
    
 \end{abstract}

\maketitle
\footnotetext[1]{N.\ Brady was partially  supported by NSF grant no.\ 
DMS-0505707}  

Hyperbolic groups contain only finitely many conjugacy classes of finite subgroups (see \cite{Bogo-Ger, Brady, Brid-Haef}).  Several other classes of groups share this property, including CAT($0$) groups \cite[Corollary~II.2.8]{Brid-Haef}, mapping class groups \cite{Bridson00}, Aut($\free$), Out($\free$) \cite{Culler84}, and arithmetic groups \cite{Borel63}.  Building on work of Grunewald and Platonov \cite{grun-plat},
Bridson \cite{Bridson00} showed that for any $n$, there is a subgroup of 
$SL(2n+2, \Z)$ of type ${\rm F}_n$ that has infinitely many conjugacy classes 
of elements of order $4$.  
In \cite{feighn-mess}, Feighn and Mess constructed finite extensions of 
$(\mathbb F_2)^n$ containing subgroups of type ${\rm F}_{n-1}$ with infinitely many conjugacy classes of elements of order $2$.  Their examples were realized as subgroups of the isometry group of $(\mathbb{H}^2)^n=\mathbb H^2 \times \dots \times \mathbb H^2$.
These examples were generalized considerably and were set in the context of 
right-angled Artin groups by Leary and Nucinkis \cite{LN}.
 
In this note we describe a model situation where the Morse theory 
argument of \cite{LN} can be applied.  It includes the right-angled Artin group setting from \cite{LN}, but it can also be used when the ambient group is 
not an Artin group.  We apply it to the case of a hyperbolic 
group in Theorem~\ref{hyperbolic} and to the case of a 
virtual direct product of hyperbolic groups in Theorem~\ref{cat}. 
In addition, we produce subgroups of mapping class groups, Aut($\free$),  and Out($\free$) with infinitely many conjugacy 
classes of finite-order elements by finding    
natural realizations of finite extensions 
of $(\mathbb F_2)^n$ in these groups (Theorems~\ref{th:mcg-example}, \ref{th:aut-example}, and \ref{th:out-example}). 
For mapping class groups this solves Problem 3.10 in \cite{Farb06}.

\section{counting conjugacy classes of finite-order elements in 
subgroups of $\Z$-extensions}

Proposition~\ref{prop:conj} below gives a basic method for producing examples of groups with infinitely many conjugacy classes of finite-order elements.  The proof is contained in that of Theorem 3 in \cite{LN}, where results are stated for right-angled Artin groups.
As background for the  model situation and 
Proposition~\ref{prop:conj}, one needs the 
notion of a Morse function on an affine cell complex, and the 
notion of non-positively curved  cell complexes. 
For details about  the first topic the reader 
can refer to \cite{BeBr}. 
A good treatment of the 
latter topic can be found in \cite{Brid-Haef}.

\medskip

\noindent
\textbf{Model Situation.} 
Let $X$ be a non-positively curved 
cell complex and let $f: X \to S^1$ be a 
circle-valued Morse function which induces 
an epimorphism of fundamental groups
$$
f_{\ast}: G = \pi_{1}(X) \to \Z = \pi_{1}(S^{1})
$$
 Let $K$ 
denote the kernel of this epimorphism.   

Now let $\sigma$ be a finite 
order isometry of $X$ which 
\begin{itemize}
\item fixes a vertex  $v \in X$, and 
\item acts without fixed points on the link $Lk(v,X)$.  
\end{itemize}
Further, assume that the map $f$ is $\sigma$-equivariant (for the trivial action of $\sigma$ on $S^1$).
The isometry $\sigma$ induces an automorphism of $G$, which we also denote by $\sigma$.  Then $f_\ast$ is $\sigma$-equivariant, and it follows that 
$\sigma$ leaves $K$ invariant.  We form the semi-direct product 
$G \rtimes \langle \sigma \rangle$, and extend $f_\ast$ to this group 
by mapping $\sigma$ trivially. 
Note that $G \rtimes \langle \sigma \rangle$ can be expressed as 
a semi-direct product with $\Z$
$$
G \rtimes \langle \sigma \rangle \; = \; (K \rtimes \langle \sigma \rangle) \rtimes \Z.
$$

\begin{prop}
\label{prop:conj} Let $\sigma$, $f$, $X$, $G$ and $K$ be as 
described in the model situation above. 
Then the group $K \rtimes \langle \sigma \rangle$ has infinitely many 
conjugacy classes of elements with the same order as $\sigma$. In fact, the 
conjugacy class of $\sigma$ in $G\rtimes \langle \sigma \rangle$ intersects 
$K \rtimes \langle \sigma \rangle$ to give infinitely many 
$K \rtimes \langle \sigma \rangle$ conjugacy classes. 
\end{prop}

\begin{proof}
Let $t\in G$ be such that $f_\ast (t)=1$.  Then $t^n \sigma t^{-n}$ is 
an element of $K \rtimes \langle \sigma \rangle$.  We show that 
$t^n \sigma t^{-n}$ is not conjugate to $t^m \sigma t^{-m}$ 
in $K \rtimes \langle \sigma \rangle$ unless 
$n = m$. 

Since $X$ is non-positively curved, its
universal cover $\widetilde X$ is a CAT($0$) space, and hence has 
unique geodesics \cite{Brid-Haef}. 
Choose a lift $\tilde \sigma : \widetilde X \to \widetilde X$ which fixes 
a vertex $x_0 \in \widetilde X$  in the preimage of the 
vertex $v \in X$ fixed by $\sigma$.  
Thus $\tilde \sigma$ is an isometry of $\widetilde X$ which fixes the 
vertex $x_{0}$. 

We argue by contradiction that this is the only point of $\widetilde X$ 
which is fixed by $\tilde \sigma$. If  $\tilde \sigma$ fixed another point  
$x_{1} \in \widetilde X$,  then it would have to fix the unique 
geodesic $[x_0x_{1}]$ (because $\tilde \sigma$ is an isometry). 
Hence $\tilde \sigma$ would fix the point of $Lk(x_0, \widetilde X)$  
determined by $[x_{0}x_{1}]$. Since 
$\tilde \sigma$ is a lift of $\sigma$, 
 this would imply that $\sigma$ 
fixes a point of $Lk(v,X)$, contradicting the hypothesis on $\sigma$. 
Thus $x_0$ is the unique fixed point of $\tilde \sigma$.

Let $\tilde{f} : \widetilde{X} \to \R$ be a lift of $f$; it is a Morse 
function on $\widetilde X$. Since $f$ is $\sigma$-equivariant 
where  $\sigma$ is defined to act trivially on $S^{1}$, 
the isometry $\tilde \sigma$ acts on $\widetilde X$ preserving 
height.

\begin{figure}[h]
\begin{center}
\psfrag{x0}{{\small $t^{m}x_{0}$}}
\psfrag{tx0}{{\small $t^{n}(x_{0})$}}
\psfrag{fx0}{{\small $\tilde{f}(x_{0})+m$}}
\psfrag{fx0n}{{\small $\tilde{f}(x_{0})+n$}}
\psfrag{0}{{\small $0$}}
\psfrag{h}{{\small $h$}}
\psfrag{f}{{\small $\tilde{f}$}}
\psfrag{x}{{\small $\widetilde{X}$}}
\psfrag{r}{{\small $\R$}}
\includegraphics[width=2.6in]{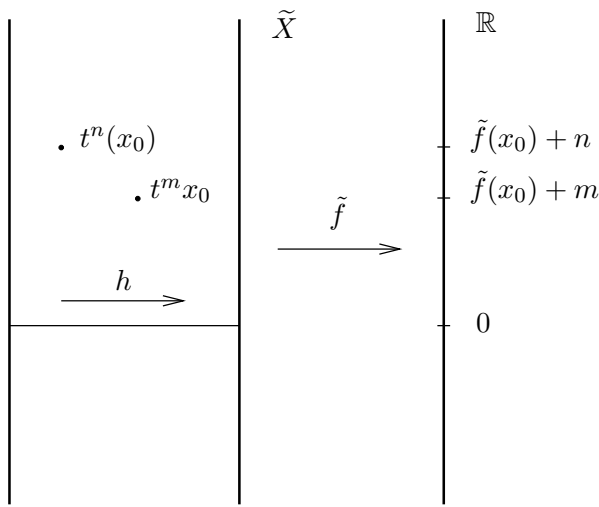}
\end{center}
\caption{Picture proof of Proposition~\ref{prop:conj}.}\label{fig:morse}
\end{figure}

Recall that we chose $t \in G$ so that $f_{\ast}(t) = 1$. The element 
$t$ acts as a deck transformation on $\widetilde X$, and so 
 $t^k\tilde \sigma t^{-k}$ 
is a lift of $\sigma$ 
with unique fixed point $t^k(x_0)$. Note that 
$$
\tilde{f}( t^k(x_0)) = \tilde{f}(x_0) + k\, .
$$  
However, each  $h \in K\rtimes \langle \sigma \rangle$ acts \emph{horizontally}    on $\widetilde X$; that is, $\tilde{f}(h(x)) = \tilde{f}(x)$ for  $x \in \widetilde{X}$. In particular $h$ cannot move the point 
$t^m x_0$ to the point $t^n x_0$ for $m\not=n$. Therefore $t^n\sigma t^{-n} \not= 
h(t^m\sigma t^{-m})h^{-1}$ for $m \not=n$, and so the elements $t^{n}\sigma t^{-n}$ each represent distinct conjugacy 
classes in $K \rtimes \langle \sigma \rangle$.  
\end{proof}

The following family of examples provides illustrations of 
Proposition~\ref{prop:conj} in all dimensions. In Section~\ref{sec:mcg} 
we shall see how to embed these examples into various classical groups.

\begin{ex} \label{ex:free}
Let $X$ be the direct product of $n$ copies of the wedge of two circles, $S^1 \vee S^1$. Then $G = \pi_1(X)$, is isomorphic to $(\mathbb F_2)^n$.
Define $f: X \to S^1$ by mapping each of the $2n$ circles homeomorphically 
around $S^1$ and extending linearly over the higher skeleta. 
Thus $f$ is a circle valued Morse function, which lifts to a 
Morse function $\tilde f$ on the universal cover $\widetilde X$. 
 The cover $\widetilde X$ is tiled by cubes. On each cube the  
map $\tilde f$ attains a maximum value at one 
vertex, and attains a minimum at the diametrically opposite 
vertex. 

Define an order $2$ automorphism $\sigma$ of $X$ as follows. 
In the case $n=1$, $\sigma$ simply interchanges the two circles. In the case 
$n \geq 2$, extend this to a diagonal action. Note that $\sigma$ 
fixes the unique vertex of $X$. The link of this vertex 
is the $n$-fold join of a set of four points.  Since $\sigma$ acts on 
a single set of four points by interchanging points in pairs, the action 
on the whole link has no fixed points. Clearly, $f$ is $\sigma$-equivariant.   

The complex $X$ is non-positively curved since it is the $n$-fold 
direct product (with the product metric) of the non-positively 
curved 1-complex $S^{1}\vee S^{1}$. If the two circles of $S^{1} 
\vee S^{1}$ are chosen to be isometric (eg.\ both are 
$\R/\Z$)  then the map 
$\sigma$ which interchanges them can be taken to be an 
isometry. Similarly, the diagonal map $\sigma$ acting 
on $X$ can be taken to be an isometry. Thus all the hypotheses 
of Proposition~\ref{prop:conj} are satisfied. 

By Proposition \ref{prop:conj}, the kernel of the map $G \rtimes \langle \sigma \rangle \to \Z$ has infinitely many conjugacy classes of elements of order $2$. .  

Note that these kernels are of type ${\rm F}_{n-1}$ but not of type ${\rm F}_n$ \cite{BeBr}. 
Examples similar to these  were given in 
Feighn-Mess \cite{feighn-mess}, and were generalized 
considerably by Leary-Nucinkis \cite{LN}.  
\end{ex}

\section{Conjugacy classes in subgroups of CAT(0) groups}
\label{sec:cat}

In this section we give two applications of Proposition~\ref{prop:conj}. 
In the first example, the 
fundamental group of the complex $X$ is hyperbolic, and the kernel of the map to $\Z$ 
is finitely generated but not finitely presented. The group 
$G = \pi_1(X) $ is extended by a 
carefully chosen finite-order automorphsim. 

\begin{theorem}\label{hyperbolic}
There exist hyperbolic groups containing 
finitely generated  subgroups which have infinitely 
many conjugacy classes of finite order elements.  
\end{theorem}

For the second example, we take the 
direct product of two copies of the complex $X$ of the first example, 
and take the diagonal action of the finite-order 
automorphism. As in the case of the direct product of free groups (in Example~\ref{ex:free} 
above) the finiteness properties of the kernel improve. In this case the 
kernel is of type ${\rm F}_3$ but not ${\rm F}_4$. The ambient CAT(0) group has $\Z^2$ 
subgroups, but not $\Z^3$ subgroups. 

\begin{theorem}\label{cat}
There exist CAT(0) groups with no $\Z^{3}$ subgroups, 
which contain 
finitely presented (in fact type F$_{3}$)   subgroups 
which have infinitely 
many conjugacy classes of finite order elements.  
\end{theorem}

\begin{rk}
It would be very  interesting to find a hyperbolic group which contains a finitely presented subgroup with infinitely many conjugacy classes
of finite-order elements. The subgroup in question could not be a hyperbolic group. 
\end{rk}
 
\medskip
\noindent
\textbf{The hyperbolic example.}
In this subsection we prove Theorem~\ref{hyperbolic} by construction. 
The construction produces one example, but it should be clear 
how to vary the construction to obtain other examples. 

\begin{figure}[h]
\begin{center}
\psfrag{i1}{{\small $i+1$}}
\psfrag{i}{{\small $i$}}
\psfrag{i2}{{\small $i+2$}}
\psfrag{i3}{{\small $i+3$}}

\includegraphics[width=2.4in]{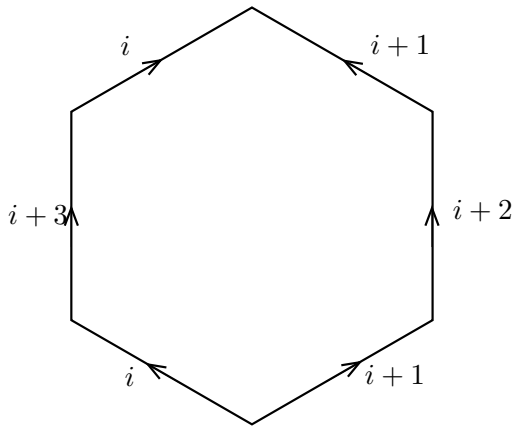}
\end{center}
\caption{A typical 2-cell of $X$.}\label{fig:2cell}
\end{figure}

Start with the group $G = \pi_{1}(X)$, where $X$ is a $2$-complex consisting of one vertex, eight $1$-cells and 
eight hexagonal $2$-cells as shown in Figure~\ref{fig:2cell}.

Define a circle-valued Morse function $f : X \to S^1$ by mapping each oriented $1$-cell homeomorphically 
around $S^1$, and extending linearly over the $2$-skeleton.  In the 
universal cover $\widetilde X$, the Morse function $\tilde f$ 
``projects'' a typical 2-cell 
(as shown in Figure~\ref{fig:2cell}) horizontally onto a segment 
of length 3 in $\R$. 

\begin{figure}[h]
\psfrag{1p}{{\small $1^+$}}
\psfrag{1m}{{\small $1^-$}}
\psfrag{2m}{{\small $2^-$}}
\psfrag{3m}{{\small $3^-$}}
\psfrag{4m}{{\small $4^-$}}
\psfrag{5m}{{\small $5^-$}}
\psfrag{6m}{{\small $6^-$}}
\psfrag{7m}{{\small $7^-$}}
\psfrag{8m}{{\small $8^-$}}
\begin{center}
\includegraphics[width=2.4in]{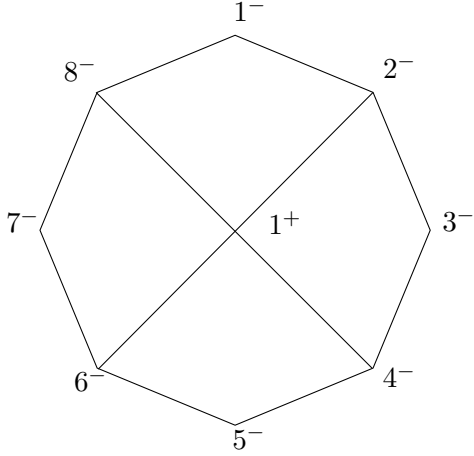}
\end{center}
\caption{Part of the link of the vertex of $X$.}
\label{fig:partiallink}
\end{figure}

The ascending  link is a circle with 8 vertices labelled $i^{-}$, 
$1\leq i \leq 8$, and 1-cells from $i^{-}$ to $(i+1)^{-}$ ($i \mod 8$).
 The descending link is described similarly, with $i^{+}$ in place of 
 $i^{-}$ above.  The full link has the following additional edges:  for each $i$, there are 
edges connecting the vertex $i^+$ to $(i\pm 1)^-$ and $(i\pm 3)^-$ (mod $8$).   Figure \ref{fig:partiallink}  shows how the vertex 
$1^{+}$ of the descending link is connected to 4 vertices of the 
ascending link.

The complex $X$ is given a piecewise hyperbolic metric by making every 
$2$-cell a regular right-angled hyperbolic hexagon.  
From Figure~\ref{fig:partiallink}, 
we see that there are no cycles of combinatorial length less 
than $4$ (or metric length less than $2\pi$).
By the large link condition, $X$ is a non-positively curved complex; the universal cover $\widetilde X$ is a CAT(-1) space, 
and so $G$ is hyperbolic. See \cite{Brid-Haef} for details about 
CAT(-1) spaces and link conditions.

Define the finite-order 
automorphism $\sigma$ of $X$ to be the cellular 
homeomorphism of the $2$-complex which fixes the vertex and cyclically permutes the oriented $1$-cells (and hence also cyclically permutes the $2$-cells). 
Note that $f$ is $\sigma$-equivariant, once we define $\sigma$ 
to act trivially on $S^1$. By construction, $\sigma$ acts without fixed points (as an order $8$ rotation) on the ascending and the descending links 
 of the single vertex of $X$. The remaining edges of the link connect 
 $+$ vertices to $-$ vertices. These are freely permuted by $\sigma$ 
 since $\sigma$ freely permutes $+$ vertices (and freely 
 permutes $-$ vertices). Finally, 
$\sigma$ can easily be taken to be an isometry of $X$. 
Thus all the hypotheses of Proposition~\ref{prop:conj} are 
satisfied. 

The fact that the ascending and descending links are each isomorphic to $S^1$ implies that the kernel $K$ of the 
map $G \to \Z$ is finitely generated but not finitely presented 
by Theorem~4.7(1) of \cite{Br} .  
Now the finite extension $G \rtimes \langle \sigma \rangle$ is virtually  $G$, and hence is also hyperbolic.  
The kernel $K\rtimes \langle \sigma \rangle$ is virtually  $K$, and so is also fintely generated 
but not finitely presented.
By Proposition~\ref{prop:conj} $K\rtimes \langle \sigma \rangle$ has infinitely many conjugacy classes of elements of order $8$.  \qed

\medskip
\noindent
\textbf{Increasing the finiteness properties of the kernel.} 
In this subsection, we give a proof of Theorem~\ref{cat} by 
construction.

Start with the hyperbolic group $G$, the 2-complex $X$, and the 
Morse function $f: X \to S^{1}$ of the 
preceding example. 

The group $G \times G$ is the fundamental group of the product space $X \times X$. Consider the composition of the product map 
$$
f \times f \, : \, X \times X \; \to \; S^1 \times S^1
$$
with the standard ``linear map'' $S^1 \times S^1 \to S^1$ (covered by the linear map $\R^{2} \to \R : (x,y) \mapsto x+y$). 
This is a circle valued Morse function on $X \times X$, whose ascending (resp.\ descending) link is simply a join of two copies of 
the ascending (resp.\ descending) link of $f$. Topologically, these ascending and 
descending links are 3-spheres. Hence the kernel $K$ of the induced map 
$G \times  G \to \Z$ is of type ${\rm F}_3$ but not of type ${\rm F}_4$ (or even ${\rm FP}_4$) by 
Theorem~4.7(1) of \cite{Br}.

Since $X$ is non-positively curved, the product complex 
$X\times X$ with the product metric is also non-positively 
curved. Hence, the universal cover $\widetilde X \times 
\widetilde X$ is a CAT(0) metric space. The group 
$G \times G$ is a CAT(0) group. Since it is the direct 
product of two hyperbolic groups, the maximal rank 
of its free abelian subgroups is 2. 

The order eight isometry $\sigma : X \to X$ defines an isometry $X\times X \to X\times X$ (acting in a diagonal fashion) which 
we also denote by $\sigma$. The Morse function on 
$X \times X$ is  $\sigma$-equivariant, once  $\sigma$ is defined 
to act trivially on $S^{1}$. 

The isometry $\sigma$ fixes the unique vertex of $X \times X$. 
We see that $\sigma$ acts without fixed points on the link of this vertex as 
follows. First note that the ascending (resp.\ descending) link 
of the Morse function on $X \times X$ is a 3-sphere, expressed 
as the join of two circles, each of combinatorial size eight. Since 
$\sigma$ acts as an order eight rotation on each circle, 
it freely permutes all the simplices of these 3-spheres. 
Finally, since  
a general simplex of  the link will be a join of a simplex of the 
ascending link and a simplex of the descending link, we 
see that $\sigma$ freely permutes all simplices of the link. 

Proposition~\ref{prop:conj} implies that the kernel 
$K \rtimes \langle \sigma \rangle$ of the induced map 
$$(G \times G) \rtimes \langle \sigma \rangle \to \Z$$ 
has infinitely many conjugacy classes of elements of order eight. 

Finally, note that the group 
$(G \times G)\rtimes \langle \sigma \rangle$ is a CAT(0) 
group, because $\sigma$ is a finite order 
isometry  of $X\times X$. 
Also, the kernel $K \rtimes \langle \sigma \rangle$ 
is virtually  $K$, and so shares the same 
finiteness properties (F$_{3}$ but not F$_{4}$).   \qed

\section{Mapping class groups, Aut($\free$), and
Out($\free$)}\label{sec:mcg}

In this section we show that the group $(\mathbb{F}_2)^n \rtimes
\langle\sigma\rangle$ from Example \ref{ex:free} can be embedded in
the mapping class group of a surface of sufficiently high genus (this
solves Problem 3.10 in \cite{Farb06}), in Aut($\free$), and
Out($\free$), where the free group has sufficiently high rank. 

 Let
$\Sigma_{g}$ be a closed orientable surface of genus $g \geq 3$ and
$\MCG(\Sigma_{g})$ the mapping class group of $\Sigma_{g}$.

\begin{theorem}\label{th:mcg-example}
There exists a subgroup of type $\mathrm{F}_{g-3}$ in $\MCG(\Sigma_{g})$ that
contains infinitely many conjugacy classes of elements of order 2.
\end{theorem}

To ensure that the example is finitely presented, we need $g \geq 5$.

\begin{figure}[ht]
\psfrag{dots}{\scriptsize $\cdots$}
\psfrag{ai}{\scriptsize $a_{i}$}
\psfrag{ai+}{\scriptsize $a_{i+1}$}
\psfrag{bi}{\scriptsize $b_{i}$}
\psfrag{bi+}{\scriptsize $b_{i+1}$}
\psfrag{s}{\scriptsize $\sigma$}
\includegraphics[width=12cm]{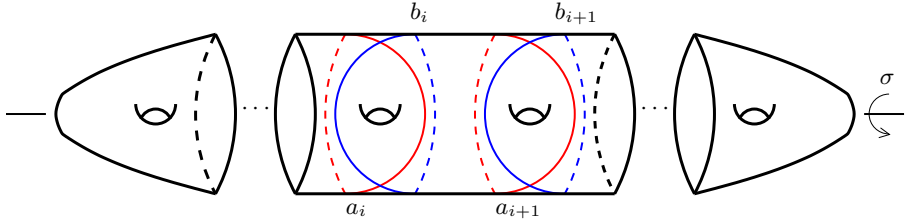}
\caption{The curves from Theorem
\ref{th:mcg-example}.}\label{fig:mcg-example}
\end{figure}

\begin{proof}
Let $a_{i},b_{i}$, $i = 1,\ldots,g-2$ be the simple closed curves
pictured in Figure \ref{fig:mcg-example} and let $T_{a_{i}},T_{b_{i}}$
denote the Dehn twists about these curves.  The only intersections
between these curves are between $a_{i}$ and $b_{i}$, which intersect
twice.  Hence by \cite[Theorem~7]{Thurston88}, the subgroup of
$\MCG(\Sigma_{g})$ generated by $T_{a_{i}},T_{b_{i}}$ is isomorphic to
$(\mathbb{F}_{2})^{g-2}$.

Let $\sigma$ be the hyperelliptic involution of $\Sigma_{g}$ such that
$\sigma(a_{i}) = b_{i}$, $\sigma(b_{i}) = a_{i}$ for all $i$.  Then
$\sigma T_{a_{i}}\sigma = T_{b_{i}}$, and hence the subgroup $H =
\langle T_{a_{i}},T_{b_{i}},\sigma \rangle < \MCG(\Sigma_{g})$ is
isomorphic to $(\mathbb{F}_{2})^{g-2} \rtimes \Z_{2}$.  The kernel
from Example \ref{ex:free} is the desired subgroup.
\end{proof}

Letting $\Sigma_{g}$ have punctures that are symmetric with respect to
the hyperelliptic involution $\sigma$, we get examples of subgroups
with infinitely many conjugacy classes of finite subgroups in
$\Out(\free)$.  The smallest finitely presented example occurs when
$g=3$ and there are $2$ punctures, this provides an example in
$\Out(\mathbb{F}_{7})$.  We give an alternative example below.

\begin{theorem}\label{th:aut-example}
There exists a subgroup of type $\mathrm{F}_{\lfloor k/2\rfloor-1}$ in
$\Aut(\free)$ that contains infinitely many conjugacy classes of
elements of order 2.
\end{theorem}

To ensure that the example is finitely presented, we need $k
\geq 6$. Note that $\lfloor k/2\rfloor$ denotes the floor of 
$k/2$.

\begin{proof}
Fix a basis $x_{1},\ldots,x_{k}$ of $\free$.  For $i =
1,\ldots,\lfloor k/2 \rfloor$ define the following automorphisms:
$$\begin{array}{crlcrl}
    \phi_{i}\co & x_{2i-1} \mapsto & x_{2i-1}^{2}x_{2i} & \psi_{i}\co 
    & x_{2i-1} \mapsto & x_{2i}x_{2i-1} \\
    & x_{2i} \mapsto & x_{2i-1}x_{2i} & & x_{2i} \mapsto & 
    x_{2i}^{2}x_{2i-1} \\
    & x_{j} \mapsto & x_{j} \mbox{ if } j \neq 2i-1,2i & & x_{j}
    \mapsto & x_{j} \mbox{ if } j \neq 2i-1,2i
\end{array}$$
For $k = 2$, the image of $\phi_{1},\psi_{1}$ in $\GL(2,\Z)$ are the
matrices ${2 \ 1 \brack 1 \ 1}$ and ${1 \ 1 \brack 1 \ 2}$.  There
exist integers $m,n\geq 1$ such that the group $\langle {2 \ 1 \brack
1 \ 1}^{m}, {1 \ 1 \brack 1 \ 2}^{n} \rangle$ is free.  Letting $N$
denote the larger of $m,n$ we see that $\langle {2 \ 1 \brack 1 \
1}^{N}, {1 \ 1 \brack 1 \ 2}^{N} \rangle$ is free.  Hence, as free
groups are Hopfian, the subgroup $\langle \phi_{1}^{N},\psi_{1}^{N}
\rangle < \Aut(\mathbb{F}_{2})$ is isomorphic to $\mathbb{F}_{2}$.  As
$\phi_{i}$ commutes with $\phi_{j},\psi_{j}$ when $j \neq i$, we see
that the group generated by the automorphisms
$\phi^{N}_{i},\psi^{N}_{i}$ is isomorphic to
$(\mathbb{F}_{2})^{\lfloor k/2 \rfloor}$.

As in the case for the mapping class group, we have an involution
$\sigma \in \Aut(\free)$ defined by $\sigma(x_{2i-1}) =
x_{2i},\sigma(x_{2i}) = x_{2i-1}$ for $i = 1,\ldots,\lfloor k/2
\rfloor$ and $\sigma(x_{k}) = x_{k}$ if $k$ is odd.  It can easily be
checked that $\sigma\phi_{i}\sigma = \psi_{i}$ and hence the subgroup
$H = \langle \phi^{N}_{i},\psi^{N}_{i},\sigma \rangle < \Aut(\free)$
is isomorphic to $(\mathbb{F}_{2})^{\lfloor k/2 \rfloor} \rtimes
\Z_{2}$.  The kernel from Example \ref{ex:free} is the desired
subgroup.
\end{proof}

\begin{theorem}\label{th:out-example}
There exists a subgroup of type $\mathrm{F}_{\lfloor k/2 \rfloor -1}$ in
$\Out(\free)$ that contains infinitely many conjugacy classes of
elements of order 2.
\end{theorem}

To ensure that the example is finitely presented, we need $k
\geq 6$.

\begin{proof}
We claim that the subgroup $H$ from Theorem \ref{th:aut-example} does
not contain any nontrivial inner automorphisms.  This implies that the
subgroup of $\Out(\free)$ generated by the images of
$\phi^{N},\psi^{N},\sigma$ is isomorphic to $H$ and hence the kernel
from Example \ref{ex:free} is the desired subgroup.

Suppose some composition $\rho$ of the
$\phi^{N}_{i},\psi^{N}_{i},\sigma$ is an inner automorphism.  Since
the image of $x_{1}$ by any automorphism in $H$ is a word in
$x_{1},x_{2}$, if $\rho(x_{1}) = xx_{1}x^{-1}$, then $x$ is a word in
$x_{1},x_{2}$.  Also, since the image of $x_{3}$ by any automorphism
is a word in $x_{3},x_{4}$, we see that if $\rho(x_{3}) =
xx_{3}x^{-1}$ then $x$ is a word in $x_{3},x_{4}$.  Therefore $x$ is
the identity and $\rho$ is trivial.
\end{proof}

\bibliographystyle{siam}
\bibliography{bibliography}

\end{document}